\documentclass[10pt]{amsart}

\usepackage{amssymb}
\usepackage{mathtools}
\usepackage{amsaddr}
\usepackage[a4paper, total={6in, 9in}]{geometry}
\usepackage{tikz}

\theoremstyle{plain}
\newtheorem{thm}{Theorem}[section]
\newtheorem{lem}[thm]{Lemma}
\newtheorem{prop}[thm]{Proposition}
\newtheorem{cor}{Corollary}
\newtheorem*{Bthm}{Birkhoff's Theorem}

\theoremstyle{definition}
\newtheorem{defn}{Definition}[section]

\theoremstyle{remark}
\newtheorem{rem}{Remark}

\newcommand{\ihat}{\breve{\imath}}
\newcommand{\jhat}{\breve{\jmath}}
\newcommand{\diag}{\mathrm{diag}}
\newcommand{\sign}{\mathrm{signum}}
\newcommand{\NtimesN}{N \times N}
\newcommand{\Rpos}{\mathbb{R}_{\geq 0}^{\NtimesN}}

\foreach \a in {q,w,e,r,t,y,u,i,o,p,a,s,d,f,g,h,j,k,l,z,x,c,v,b,n,m,Q,W,E,R,T,Y,U,I,O,P,A,S,D,F,G,H,J,K,L,Z,X,C,V,B,N,M}{
  \expandafter\xdef\csname v\a\endcsname {
    {\noexpand\mathbf{\a}}
  }
}

\title[From Local Updates to Global Balance]{From Local Updates to Global Balance: A Framework for Distributed Matrix Scaling}
\thanks{G.~Aletti and G.~Naldi are members of the Italian Group ``Gruppo Nazionale per il Calcolo Scientifico'' of the Italian Institute ``Istituto Nazionale di Alta Matematica''.}

\author[G.~Aletti and G.~Naldi]{Giacomo Aletti and Giovanni Naldi}
\address{DESIRE center, Università degli Studi di Milano, Milan, Italy}
\email{giacomo.aletti@unimi.it, giovanni.naldi@unimi.it}
\thanks{G.~Naldi is the corresponding author.}


\date{\today}

\dedicatory{}

\begin{document}

\begin{abstract}
This paper investigates matrix scaling processes in the context of local normalization algorithms and their convergence behavior. Starting from the classical Sinkhorn algorithm, the authors introduce a generalization where only a single row or column is normalized at each step, without restrictions on the update order. They extend Birkhoff's theorem to characterize the convergence properties of these algorithms, especially when the normalization sequence is arbitrary. A novel application is explored in the form of a Decentralized Random Walk (DRW) on directed graphs, where agents modify edge weights locally without global knowledge or memory. The paper shows that such local updates lead to convergence towards a doubly stochastic matrix, ensuring a uniform stationary distribution across graph vertices. These results not only deepen the theoretical understanding of matrix scaling but also open avenues for distributed and agent-based models in networks.
\end{abstract}

\maketitle


\section{Introduction}

The matrix scaling problem involves the transformation of a given nonnegative matrix using diagonal matrices so that the resulting matrix has prescribed row and column sums \cite{sinkhorn1964a}. More formally, given a nonnegative matrix $A \in \mathbb{R}^{m \times n}$ and vectors $r \in \mathbb{R}^m$ and $c \in \mathbb{R}^n$ of target row and column sums, the goal is to find diagonal matrices $D_1$ and $D_2$ such that $B = D_1 A D_2$ satisfies the marginal constraints.

Many researchers have been aware of parallel efforts on matrix scaling, yet few have attempted to provide comprehensive overviews. This highlights the complexity of the field’s history, which often makes it difficult to determine whether a particular problem has already been addressed. Notable efforts to document the historical development include the recent works of \cite{kalantari2008a,pukelsheim2009a,Idel_review16}. 

Today, matrix scaling remains a well-studied problem, central to both theoretical exploration and practical applications in various disciplines. In statistics, matrix scaling appears in contingency table adjustments and maximum likelihood estimation under marginal constraints \cite{ireland1968a,jaynes1957a}. In economics, it underpins input-output models to balance supply and demand across sectors. In transportation planning, it is used to estimate trip distributions using known marginals. Furthermore, in numerical analysis, matrix scaling is often applied to improve the numerical stability of algorithms.

In recent years, matrix scaling has also been influential in computer science and quantum information theory \cite{ivanyos2015a}. It plays a role in operator scaling and positive map transformations \cite{gurvits2004a}, with implications for complexity theory and quantum entropy inequalities. The wide applicability of this problem makes it a cornerstone topic across multiple disciplines.

% \subsection*{Historical Development of the Problem}

The origins of matrix scaling are multifaceted, stemming from distinct and largely unconnected fields. The earliest documented use dates to Kruithof (1937), who devised a method to adjust telecommunication matrices based on marginal forecasts. This method, now known as the Kruithof projection \cite{kruithof1937a}, laid the groundwork for iterative approaches to matrix balancing.

In 1940, Deming and Stephan \cite{deming1940a} introduced the iterative proportional fitting procedure (IPFP) to align observed contingency tables with known marginals \cite{stephan1942a}. Though their method was widely adopted, their justification was ad hoc and flawed, as later corrected by Stephan himself. The terminology ``RAS method'' became prevalent in economics due to Richard Stone’s contributions in the early 1960s, where it was used to update input-output tables by balancing row and column sums \cite{stone1962a}.

Another foundational development occurred in statistical physics with Schrödinger’s 1931 inquiry into Brownian motion \cite{schroedinger1931a}, which can be interpreted as a continuous version of matrix scaling. Fortet (1940) addressed this problem with fixed-point techniques \cite{fortet1940a}, making early strides in establishing convergence, albeit in challenging notation.

The formal mathematical framing of the problem took shape in the 1960s through the work of Sinkhorn, who proved that any square matrix with strictly positive entries can be scaled to a doubly stochastic matrix via diagonal equivalence \cite{sinkhorn1964a}. His algorithm, along with contributions by Knopp, demonstrated convergence properties and inspired subsequent research \cite{sinkhorn1967b}.

Over the ensuing decades, research proliferated. Macgill (1977) surveyed computational methods \cite{macgill1977a}, while Rothblum and Schneider (1989) collected hundreds of related publications \cite{rothblum1989b}. Reviews by Fienberg (1970), Schneider and Zenios (1990), and Kalantari and Khachiyan (1996) have traced the development of the problem across various domains, revealing both the depth and breadth of the matrix scaling literature \cite{fienberg1970a,schneider1990b,kalantari1996b}.

% \section*{3. Summary of Approaches to Matrix Scaling}

Matrix scaling has inspired a rich set of mathematical tools for proving the existence and constructing transformations such as $D_1AD_2$ to specified row and column sums. These methods, many interlinked, provide a comprehensive landscape of techniques.

%\subsection*{Logarithmic Barrier Function}

One key method is the \textbf{logarithmic barrier function}, initially proposed by \cite{gorman1963a} and rigorously discussed by \cite{kalantari1996b}. The function
\[
g(x, y) = y^T A x - \sum_i c_i \ln x_i - \sum_i r_i \ln y_i
\]
has a stationary point if and only if $A$ can be scaled. \cite{marshall1968a} later framed this in Lagrangian terms, linking scaling to convex optimization.

% \subsection*{Nonlinear Perron-Frobenius Theory}

Another method uses \textbf{nonlinear Perron-Frobenius theory}, where scaling is reinterpreted as finding fixed points of the Menon operator $T(x) = c / (A^T (r / Ax))$, originally studied by \cite{brualdi1966a}. The operator is contractive under Hilbert’s metric, and Brouwer’s theorem assures fixed-point existence in compact domains. The work of Brualdi, Parter, and Schneider laid the foundations for this theory, connecting it with eigenvalue problems and dynamical systems \cite{brualdi1966a}.

\textbf{Entropy minimization} interprets matrix scaling as minimizing the Kullback-Leibler divergence between the scaled matrix and the original matrix under marginal constraints. Csiszár (1975) and Franklin and Lorenz (1989) established strong theoretical foundations using the framework of information geometry \cite{csisz1975a,franklin1989a}. These methods reveal the RAS algorithm as a coordinate ascent method in a dual space.

\textbf{Convex programming} builds on this by converting the optimization problem using exponential substitutions. \cite{kalantari1996a} show this results in a strictly convex function $f(x, h)$ whose minimum yields the scaling. Duality reveals this as the Wolfe or Lagrangian dual of the entropy formulation \cite{macgill1977a, balakrishnan2004a}.

\textbf{Topological approaches}, pioneered by Bapat (1982) and later expanded by Raghavan (1984) and Friedland (2016) \cite{bapat1982a,raghavan1984a,friedland2016a}, focus on proving the existence of solutions using fixed-point theorems such as Brouwer and Kakutani. These methods emphasize the continuity and surjectivity of scaling transformations, independent of algorithmic constructs. Degree theory and homotopy arguments ensure the robustness of these existence results.

Letac (1974) proposed a general \textbf{nonlinear projection} theorem that embeds matrix scaling into the broader context of convex projection and functional analysis \cite{letac1974a}. This approach is highly abstract but unifies several known results by showing that solutions correspond to projections onto convex sets under exponential mappings. It has deep implications for optimization theory and mathematical economics.

\smallskip

In the general RAS problem, the goal is to find diagonal matrices $D_1$ and $D_2$ such that $B = D_1 A D_2$ has prescribed row sums $r$ and column sums $c$. The doubly stochastic case is a special instance where $r = c = \mathbf{1}$, i.e., all row and column sums equal 1. This symmetry simplifies the analysis and makes it possible to apply strong theoretical results like Sinkhorn’s theorem \cite[Theorem]{Sinkhorn67}.
Moreover, every positive matrix that is scalable can be transformed into a doubly stochastic matrix. This allows many core results (e.g., uniqueness, convergence) to be proven in the doubly stochastic setting and generalized to the broader case.
In addition, in Markov chains and network diffusion, doubly stochastic matrices correspond to random walks with uniform stationary distribution, a desirable feature in decentralized systems and fairness-based models. 

Summing up, we focus on the doubly stochastic case that provides cleaner mathematical properties and robust convergence guarantees. Once mastered, we extends our framework naturally to general marginals in a optimal transport application and in network theory.
From now on, given a nonnegative matrix $A$, the scaling objective is to find diagonal matrices $D_1$ and $D_2$ such that $B = D_1 A D_2$, 
$B \mathbf{1} = \mathbf{1}$, and  $\mathbf{1}^\top B = \mathbf{1}^\top$.

A foundational result in matrix theory that plays a pivotal role is the Birkhoff’s theorem \cite{birkhoff1957a}. The theorem states that every $n \times n$ doubly stochastic matrix can be represented as a convex combination of permutation matrices. Mathematically, if $B$ is a doubly stochastic matrix, then there exist permutation matrices $P_1, P_2, \ldots, P_k$ and nonnegative weights $\lambda_1, \lambda_2, \ldots, \lambda_k$ summing to 1 such that:
\[
B = \sum_{i=1}^k \lambda_i P_i.
\]

This characterization has several important implications for matrix scaling. First, it provides a geometric interpretation: the set of doubly stochastic matrices forms a convex polytope, known as the Birkhoff polytope, whose vertices are precisely the permutation matrices. This geometric structure is crucial when considering optimization over the space of doubly stochastic matrices, as many scaling algorithms aim to project or approximate within this convex set.

Starting with a matrix $A$, the alternating Sinkhorn algorithm alternately normalizes the rows and columns: at each iteration, it first scales each row so the row sums equal one, then scales each column to achieve column sums of one \cite{sinkhorn1967b}. Although the existence of two diagonal matrices $D_1$ and $D_2$ such that $B = D_1 A D_2$ is closely linked to the pattern of $A$, it is well known that Sinkhorn's algorithm may still converge to a matrix $B$ even when an exact scaling solution does not exist. One of the contributions of this paper is a complete characterization of the class of limit matrices to which Sinkhorn's algorithm may converge. In particular, we identify the submatrix $A^*$ of $A$ such that the limiting matrix $B$ satisfies the scaling relation $B = D_1 A^* D_2$.

The main theoretical result of this paper is an extension of existing results to a broader class of local normalization algorithms. In these algorithms, only a single row or column is normalized at each step, with no constraints on the order or method of selection. The final characterization is achieved through an appropriate generalization of Birkhoff’s theorem.
This class of local algorithms is novel in the literature, and the results offer a characterization of their convergence in the context of distributed and agent-based normalization processes.  

A significant initial application of these new results is in the study of random walks on directed, strongly connected graphs, where local adjustments can lead to a doubly stochastic limiting matrix. In this context, it is important to recall that a uniform invariant distribution is a defining property of doubly stochastic matrices. Since the frequency of visits to each vertex in a random walk is proportional to its eigenvector centrality score, this can result in polarization. We show that a single agent, by exploring the graph and locally modifying the transition probabilities of both outgoing and incoming edges at the vertices it visits --without retaining memory of the path-- can asymptotically achieve a uniform distribution of time spent across all vertices. This result generalizes in a straightforward manner to a finite number of agents.

The structure of the paper is as follows.
The second section introduces the mathematical framework based on directed graphs with weighted edges, where the adjacency matrix can be normalized to become row, column, or doubly stochastic. Convex sets and distance measures are defined to quantify how close a matrix is to being doubly stochastic. It is shown that local row or column normalization consistently reduces this distance, providing the foundation for the convergence analysis.
In the third section, we generalize the Sinkhorn process by allowing sequences of arbitrary local updates. We prove convergence under mild assumptions using an extended version of the Birkhoff-von Neumann theorem. This extension accounts for partial supports and submatrices and provides a characterization of the uniqueness and properties of the limiting matrix.
The fourth section explores two applications. In the first one, a single agent navigates a directed graph, locally adjusting edge weights at visited vertices—resulting in what the authors call a Decentralized Random Walk (DRW). Despite lacking global information or memory, the agent's local updates lead the system to asymptotically converge to a uniform stationary distribution, ensuring balanced exploration across all vertices.
The second application focuses on entropically regularized \textit{Optimal Transport (OT)} problems, solved efficiently via the Sinkhorn algorithm. While this approach enables fast computation, the results reveal a fundamental limitation: as the regularization parameter \( \varepsilon \to 0 \), numerical underflow can lead to the loss of support in the scaled cost matrix \( \xi = e^{-C/\varepsilon} \), preventing convergence.

%%%%%%%%%%%%%%%%%%%%%%%%%%%%%%%%%%%%%%%%%%%%%%%%%%
%% results
%%%%%%%%%%%%%%%%%%%%%%%%%%%%%%%%%%%%%%%%%%%%%%%%%%

\section{Graphs and row-column stochasticity}

Bistochastic (or doubly stochastic) matrices play a key role in modeling random walks on directed graphs. When used as the transition matrix of a random walk, a bistochastic matrix ensures that the uniform distribution is invariant: the probability of being at any given vertex remains constant over time. This property makes bistochastic matrices particularly relevant for modeling \emph{unbiased exploration} in networks. In the context of directed graphs, where edges have orientation and possibly nonuniform weights, enforcing bistochasticity through local or global normalization guarantees that the random walk does not favor certain nodes asymmetrically, leading instead to an equal long-term visitation rate across all vertices.

\begin{defn}[Directed Graph]
A \emph{directed graph} (or \emph{digraph}) is a pair \( G = (V, E) \) where \( V \) is a set of vertices and \( E \) is a set of ordered pairs of vertices, called edges. Each edge is directed from one vertex to another.\\
% A \emph{directed non-negative weighted graph} is a directed graph where each edge has a non-negative weight associated with it. Formally, it is a triplet \( G = (V, E, w) \), where \( V \) is a set of vertices, \( E \) is a set of ordered pairs of vertices (edges), and \( w: E \to \mathbb{R}_{\geq 0} \) is a weight function that assigns a non-negative real number to each edge. 
    A \emph{function $w$ on the edges of a (directed) graph} $G = (V,E)$ is a function \( w: E \to \mathbb{R} \) that assigns a real number to each edge.
    A \emph{weight function on the edges} is a non-negative edge function \( w: E \to \mathbb{R}_{\geq 0} \).\\
    A \emph{(directed) weighted graph} is a directed graph with a weight function on the edges. Formally, it is a triplet \(  (V, E, w) \), where \( (E,V) \) is a directed graph, and \( w \) is a weight function on the edges.
\end{defn}

\begin{defn}
    The \emph{matrix representation} of a directed non-negative weighted graph \( G = (V, E, w) \) is a matrix \( W \in \mathbb{R}^{\NtimesN}\), where \( n \) is the number of vertices. The entry \( W_{ij} \) represents the weight of the edge from vertex \( i \) to vertex \( j \). If there is no edge from \( i \) to \( j \), \( W_{ij} = 0 \).\\
    A \emph{row (resp.\ column) stochastic matrix} is a square matrix in which each row (resp.\ column) sums to 1, and all entries are non-negative. Formally, a matrix \( W \) is row (column) stochastic if \( \sum_{j} W_{ij} = 1 \) (resp.\ \( \sum_{j} W_{ji} = 1 \)) for all \( i \), and \( W_{ij} \geq 0 \) for all \( i, j \). It represents a non-negative weighted graph for which \( \sum_{v \in V \colon (v_1,v)\in E} w(e) = 1 \) (resp.\ \( \sum_{v \in V \colon (v,v_1)\in E} w(e) = 1 \)) for all \( v_1 \in V \). A \emph{bistochastic matrix} is a square matrix which is row and column stochastic. We denote by $\mathcal{R}$ and $\mathcal{C}$ the set of the row and column stochastic matrices, and by $\mathcal{B}= \mathcal{R} \cap \mathcal{C}$ the set of bistochastic matrices.
\end{defn}
\begin{rem}
$\mathcal{R}$, $\mathcal{C}$ and $\mathcal{B}$ are convex compact subsets of the convex cone of the non-negative matrices. 
\end{rem}

The \( L_{1,1} \) norm of a matrix, also known as the entrywise \( L_1 \) norm, is the sum of the absolute values of all the entries in the matrix. 
For a matrix \( W \) with elements \( W_{ij} \), the \( L_{1,1} \) norm is defined as:
\[
\|W\|_{1,1} = \sum_{i} \sum_{j} |W_{ij}| = \sum_{i,j} |W_{ij}|.
\]
From now on, the set of nonnegative matrix without null rows and columns will be denoted by $\Rpos$.
\begin{lem}\label{lem:d_Rd_C_repr}
    Let $W \in \Rpos $ be a non-negative matrix. The $L_1$ distance $d_{\mathcal{R}}$ of $W$ from $\mathcal{R}$ (resp.\ $d_{\mathcal{C}}$ distance from $\mathcal{C}$) is
\[
d_{\mathcal{R}} (W) = \inf_{R\in \mathcal{R}} \|W-R\|_{1,1} = \sum_{i} \Big| \sum_{j} W_{ij} -1 \Big|,
\qquad
d_{\mathcal{C}} (W) = \inf_{C\in \mathcal{C}} \|W-C\|_{1,1} = \sum_{j} \Big| \sum_{i} W_{ij} -1 \Big|,
\]
\end{lem}
\begin{rem}\label{rem:d_Rd_C_top}
From Lemma~\ref{lem:d_Rd_C_repr}, \(d_{\mathcal{R}} (W) = d_{\mathcal{C}} (W^\top)\).
\end{rem}
\begin{proof}
We prove for $d_{\mathcal{R}}$, since $C^\top\in\mathcal{R}$ whenever $C\in\mathcal{C}$. Note that, for any $R\in\mathcal{R}$,
\[
\|W-R\|_{1,1} = \sum_{i} \sum_{j} | W_{ij} -R_{ij} | 
\geq \sum_{i} \Big| \sum_{j} W_{ij} - \sum_{j} R_{ij} \Big| = \sum_{i} \Big| \sum_{j} W_{ij} -1 \Big|.
\]
Now, let $ \bar{R}_{ij} = W_{ij} / (\sum_{k} W_{ik} )$ (if $ \sum_{k} W_{ik} = 0$ then $\bar{R}_{ij} = 1/N$ for any $j$). Then
\[
\|W-\bar{R}\|_{1,1} = \sum_{i} \sum_{j} |W_{ij}-\bar{R}_{ij}| = \sum_{i} \sum_{j} \Big| \sum_{k} W_{ik} - 1 \Big| |\bar{R}_{ij}|
= \sum_{i} \Big| \sum_{k} W_{ik} - 1 \Big| \sum_{j} |\bar{R}_{ij}| = \sum_{i} \Big| \sum_{k} W_{ik} - 1 \Big| . \qedhere
\]
\end{proof}

\subsection{Row-column normalization}
\begin{defn}\label{defn:RC(W)}
The two operators $T^{(r)}_{\ihat},T^{(c)}_{\jhat}:\Rpos \to \Rpos $ ($\ihat,\jhat\in \{1,\ldots,N\}$) that normalize the $\ihat$-th row and the $\jhat$-th column are defined as
\[
T^{(r)}_{\ihat} (W) = 
\begin{cases}
    W_{ij} & \text{if $i \neq \ihat$}
    \\
    \frac{W_{ij}}{\sum_k W_{ik}} & \text{if $i = \ihat$},
\end{cases}
    \qquad
T^{(c)}_{\jhat} (W) = 
\begin{cases}
    W_{ij} & \text{if $j \neq \jhat$}
    \\
    \frac{W_{ij}}{\sum_k W_{kj}} & \text{if $j = \jhat$}.
\end{cases}
\]
For $\ihat = \jhat = 0$, define $T^{(r)}_{0}=T^{(c)}_{0} $ as the identity operator.
\end{defn}
\begin{rem}\label{rem:r_ic_j_top}
    Note that \(T^{(c)}_{\jhat} (W) = (T^{(r)}_{\jhat} (W^\top))^\top\). 
\end{rem}
When we refer a general normalizing operator, for $k \in \{r,c\}$ and $l\in\{1,\ldots,N\}$, we use $T^{(k)}_l$ to denote both $T^{(r)}_{\ihat}$ and $T^{(c)}_{\jhat}$. 
For $W\in\Rpos$, we define
\[
W^{(r)}_i = \sum_j W_{ij} - 1, \qquad W^{(c)}_j = \sum_i W_{ij} - 1,
\]
so that, by Lemma~\ref{lem:d_Rd_C_repr}
\begin{equation}\label{eq:dRdCdB}
d_{\mathcal{R}} (W) = \sum_i |W^{(r)}_i|, \qquad
d_{\mathcal{C}} (W) = \sum_j |W^{(c)}_i|.
\end{equation}
Given $f:\mathbb{N}\to \mathbb{R}$, for any $t_0,t_1\in\mathbb{N}$,
\[
\Delta_{t_0,t_1} f(\cdot) = f(t_1)-f(t_0).
\]
Now we underline the effect of the operators $T^{(k)}_l$ on $d_{\mathcal{R}}$, on $d_{\mathcal{C}}$, and on their sum $ d_{\mathcal{B}} = d_{\mathcal{R}} + d_{\mathcal{C}}$.
    Fixed $W\in \Rpos$, define $W_0= W$ and $W_1 = T^{(k)}_l(W)$. Then
\[
\begin{aligned}
\Delta_{0,1} (W_{\cdot})^{(r)}_i
= 
(W_1)^{(r)}_i - (W_0)^{(r)}_i 
& = 
\begin{cases}
\begin{cases}
    - (W_0)^{(r)}_i 
    & \text{if }i = l;
    \\
    0 & \text{if }i \neq l;
\end{cases}
& \text{if } k = r;
\\
\Delta_{0,1} (W_{\cdot})^{(r)}_i
& \text{if } k = c;
\end{cases}
\\
\Delta_{0,1} (W_{\cdot})^{(c)}_j
= 
(W_1)^{(c)}_j - (W_0)^{(c)}_j 
& = 
\begin{cases}
\Delta_{0,1} (W_{\cdot})^{(c)}_j
& \text{if } k = r;
\\
\begin{cases}
    - (W_0)^{(c)}_j 
    & \text{if }j = l;
    \\
    0 & \text{if }j \neq l;
\end{cases}
& \text{if } k = c.
\end{cases}
\end{aligned}
\]
Since the sum of row sums and column sums coincides, for $T^{(k)}_l$ we have
\begin{equation}\label{eq:DeltaW^k_l}
\begin{aligned}
    & \left.
\begin{aligned}
    &
    (W_1)^{(r)}_l = 0 ,  \qquad  \Delta_{0,1} (W_{\cdot})^{(r)}_i = 0, \ \forall i \neq l, 
    \\
    &
    (W_0)^{(r)}_l = - \Delta_{0,1} (W_{\cdot})^{(r)}_l = - \sum_j \Delta_{0,1} (W_{\cdot})^{(c)}_j , \\
    & 
    |(W_0)^{(r)}_l| = \sum_j | \Delta_{0,1} (W_{\cdot})^{(c)}_j |
    \\
    & |(W_0)^{(r)}_l| = - \Delta_{0,1} | (W_{\cdot})^{(r)}_l| = - \sum_i \Delta_{0,1} | (W_{\cdot})^{(r)}_i |     
\end{aligned}
    \right\}
    && \text{if }k = r;
    \\
    & \left.
\begin{aligned}
    & (W_1)^{(c)}_l = 0, \qquad     \Delta_{0,1} (W_{\cdot})^{(c)}_j = 0, \ \forall j \neq l, 
    \\
    & (W_0)^{(c)}_l = - \Delta_{0,1} (W_{\cdot})^{(c)}_l = - \sum_i \Delta_{0,1} (W_{\cdot})^{(r)}_i , \\
    & 
    |(W_0)^{(c)}_l| = \sum_i | \Delta_{0,1} (W_{\cdot})^{(r)}_i |
    \\
    & |(W_0)^{(c)}_l| = - \Delta_{0,1} | (W_{\cdot})^{(c)}_l| = - \sum_j \Delta_{0,1} | (W_{\cdot})^{(c)}_j | 
\end{aligned}
    \right\}
    && \text{if }k = c.
\end{aligned}
\end{equation}
Each normalization scales either enlarging or reducing the matrix terms, and hence
\begin{equation}
\begin{aligned}
& 
\big(\Delta_{0,1} (W_\cdot)^{(r)}_i\big) 
\big(\Delta_{0,1} (W_\cdot)^{(c)}_j\big) \geq 0;
&& \forall (i,j) \in \NtimesN, 
\\
&
\sum_i | \Delta_{0,1} (W_{\cdot})^{(r)}_i | 
= \Big | \sum_i \Delta_{0,1} (W_{\cdot})^{(r)}_i \Big|;
\\
&
\sum_j | \Delta_{0,1} (W_{\cdot})^{(c)}_j | = \Big| \sum_j \Delta_{0,1} (W_{\cdot})^{(c)}_j \Big|.
\end{aligned}
\end{equation}
To check the first of the two relations above, note that by \eqref{eq:DeltaW^k_l} we obtain
\[
\sum_i | \Delta_{0,1} (W_{\cdot})^{(r)}_i | =
\begin{cases}
    |0 - (W_{0})^{(r)}_l| + \sum_{i\neq l} 0 = |(W_{0})^{(r)}_l| =  \Big | \sum_i \Delta_{0,1} (W_{\cdot})^{(r)}_i \Big|
    & \text{if }k = r,
    \\
    |(W_0)^{(c)}_l|  = \Big| - \sum_i \Delta_{0,1} (W_{\cdot})^{(r)}_i  \Big|  =  \Big | \sum_i \Delta_{0,1} (W_{\cdot})^{(c)}_i \Big|
    & \text{if }k = c.
\end{cases}
\]
%    Fixed $W\in \Rpos$, define $W_0= W$ and $W_1 = T^{(k)}_l(W)$. Then, 
Summing up, by \eqref{eq:dRdCdB} and \eqref{eq:DeltaW^k_l},
\[
\begin{aligned}
& d_{\mathcal{R}} (W_{1}) - d_{\mathcal{R}} (W_{0}) =
\sum_i \Delta_{0,1} | (W_{\cdot})^{(r)}_i | = 
\begin{cases}
    -|(W_0)^{(r)}_l| & \text{if } k = r
    \\
    \sum_i \Delta_{0,1} | (W_{\cdot})^{(r)}_i | 
    & \text{if } k = c
\end{cases}
\\
& d_{\mathcal{C}} (W_{1}) - d_{\mathcal{C}} (W_{0}) =
\sum_j \Delta_{0,1} | (W_{\cdot})^{(c)}_j | = 
\begin{cases}
    \sum_j \Delta_{0,1} | (W_{\cdot})^{(c)}_j | 
    & \text{if } k = r
    \\
    -|(W_0)^{(c)}_l| 
    & \text{if } k = c
\end{cases}
\end{aligned}
\]
Recall that $ d_{\mathcal{B}} = d_{\mathcal{R}} + d_{\mathcal{C}}$; by \eqref{eq:DeltaW^k_l} we have
\begin{equation}\label{eq:deltaD_B}
d_{\mathcal{B}} (W_{1}) - d_{\mathcal{B}} (W_{0}) =
\begin{cases}
    \sum_j \Delta_{0,1} | (W_{\cdot})^{(c)}_j | -  \sum_j | \Delta_{0,1} (W_{\cdot})^{(c)}_j | 
    & \text{if } k = r;
    \\
    \sum_i \Delta_{0,1} | (W_{\cdot})^{(r)}_i | -  \sum_i | \Delta_{0,1} (W_{\cdot})^{(r)}_i | 
    & \text{if } k = c.
\end{cases}
\end{equation}
Finally, note the matrix representation
\begin{equation}\label{eq:matrTrepres}
W_1 = T^{(k)}_l(W_0) =
\begin{cases}
\begin{psmallmatrix}
1 & 0 & \cdots & 0 & \cdots & 0 & 0
\\
0 & 1 & \cdots & 0 & \cdots & 0 & 0
\\
0 & 0 & \cdot & 0 & \cdots & 0 & 0
\\
0 & 0 & \cdots & (1+(W_0)^{(r)}_l)^{-1} & \cdots & 0 & 0
\\
0 & 0 & \cdots & 0 & \cdot & 0 & 0
\\
0 & 0 & \cdots & 0 & \cdots & 1 & 0
\\
0 & 0 & \cdots & 0 & \cdots & 0 & 1
\end{psmallmatrix}
W_0 = 
    \diag \big( 1,\ldots, 1, \tfrac{1}{1+(W_0)^{(r)}_l},1,\ldots, ,1 \big) W_0
    & \text{if } k = r
    \\
W_0
\begin{psmallmatrix}
1 & 0 & \cdots & 0 & \cdots & 0 & 0
\\
0 & 1 & \cdots & 0 & \cdots & 0 & 0
\\
0 & 0 & \cdot & 0 & \cdots & 0 & 0
\\
0 & 0 & \cdots & (1+(W_0)^{(c)}_l)^{-1} & \cdots & 0 & 0
\\
0 & 0 & \cdots & 0 & \cdot & 0 & 0
\\
0 & 0 & \cdots & 0 & \cdots & 1 & 0
\\
0 & 0 & \cdots & 0 & \cdots & 0 & 1
\end{psmallmatrix}
= 
    W_0 \diag \big( 1,\ldots, 1, \tfrac{1}{1+(W_0)^{(c)}_l},1,\ldots, ,1 \big) 
    & \text{if } k = c
\end{cases}
\end{equation}

\subsection{Sequences of normalizations}
Now we deal with a sequence $(k_t,l_t)_t$ of row/column normalizations. 
The Sinkhorn–Knopp algorithm is an interesting example of such sequences with 
\[
\begin{aligned}
k_t & = 
\begin{cases}
    r & \text {if $t \in \{2kN+1, \ldots, (2k+1)N \}$, for some $k\in \mathbb{N}$};
    \\
    c & \text {if $t \in \{(2k+1)N+1, \ldots, (2k+2)N \}$, for some $k\in \mathbb{N}$};
\end{cases}
\\
l_t & = t\mod{N}, \text{ with the obvious notation }N \equiv 0.
\end{aligned}
\]

\begin{defn}
For any $t\in\mathbb{N}$, let $k_t\in\{r,c\}$, $l_t\in\{1,\ldots,N\}$ and $(T^{(k_t)}_{l_t})_{t\in\mathbb{N}}$ the corresponding sequence of normalizing operators. Given $W_0\in\Rpos$,  the sequence \(W_t = T^{(k_t)}_{l_t}(W_{t-1})\) inductively defined is called \emph{normalized sequence} generated by $(T^{(k_t)}_{l_t})_{t}$.
\end{defn}
By \eqref{eq:matrTrepres},
\begin{equation}\label{eq:matreDtRepr}
    W_t = D^{(r)}_t W_0 D^{(c)}_t ,
\end{equation}
where $D^{(r)}_t$ and $D^{(c)}_t$ are diagonal matrices with elements
\[
(D^{(r)}_t)_{ii} = \prod_{s \leq t\colon k_s = r, l_s = i} \frac{1}{1+(W_s)^{(r)}_i}, \qquad 
(D^{(c)}_t)_{jj} = \prod_{s \leq t\colon k_s = c, l_s = j} \frac{1}{1+(W_s)^{(c)}_j}.
\]
\begin{prop}\label{prop:boundNormsGen}
    Let $(W_t)_t$ be a normalized sequence generated by $(T^{(k_t)}_{l_t})_{t}$.     
    Then for any $t_0\leq t_1\in\mathbb{N}$, \(\Delta_{t_0,t_1} d_{\mathcal{B}}(W_\cdot)\leq 0\), and, in particular, $d_{\mathcal{B}}(W_t)\leq d_{\mathcal{B}}(W_0)$ for any $t$.
    In addition, for any $T\in\mathbb{N}$, we have
    \[
    \Big|
    \sum_{t=1}^T (W_t)^{(k_t)}_{l_t}
    \Big|
    \leq 2 d_{\mathcal{B}}(W_{0}), 
    \qquad
    \text{and hence}
    \qquad
    \prod_{t=1}^T \frac{1}{1+(W_t)^{(k_t)}_{l_t}} \geq \exp(-2 d_{\mathcal{B}}(W_{0})).
    \]
\end{prop}
\begin{proof}
    Since $|a|-|b| - |a-b| \leq 0$, then \eqref{eq:deltaD_B} implies \(d_{\mathcal{B}}(W_t)\leq d_{\mathcal{B}}(W_{t-1})\), for any $t$. Then \(\Delta_{t_0,t_1} d_{\mathcal{B}}(W_\cdot)\leq 0\), and, in particular, \(d_{\mathcal{B}}(W_t)-d_{\mathcal{B}}(W_0) = \Delta_{0,t} d_{\mathcal{B}}(W_\cdot)\leq 0\).

Now, by \eqref{eq:DeltaW^k_l}, note that for any $t$, since the row sums equalize the column sums,
\begin{multline*}
(W_t)^{(k_t)}_{l_t} 
= - \big( (W_t)^{(k_t)}_{l_t} - (W_{t-1})^{(k_t)}_{l_t} \big) 
= - \big( (W_t)^{(k_t)}_{l_t} - (W_{t-1})^{(k_t)}_{l_t} \big) 
- \sum_{i\neq l_t}\big( (W_t)^{(k_t)}_{i} - (W_{t-1})^{(k_t)}_{i} \big)
\\
= - \sum_i \big( (W_t)^{(k_t)}_{i} - (W_{t-1})^{(k_t)}_{i} \big)
= - \sum_i \big( (W_t)^{(r)}_{i} - (W_{t-1})^{(r)}_{i} \big) = - \sum_i \Delta_{t-1,t} (W_{\cdot})^{(r)}_{i}
\end{multline*}
and hence
\begin{multline*}
    \Big|
    - \sum_{t=1}^T (W_t)^{(k_t)}_{l_t}
    \Big|
=
    \Big|
    \sum_i \sum_{t=1}^T \Delta_{t-1,t} (W_{\cdot})^{(r)}_{i}
    \Big|
=
    \Big|
    \sum_i \Delta_{0,T} (W_{\cdot})^{(r)}_{i}
    \Big|
=
    \Big|
    \sum_i \big( (W_{T})^{(r)}_{i} -
     (W_{0})^{(r)}_{i} \big) 
    \Big|
    \\
\leq 
    \sum_i | (W_{T})^{(r)}_{i} | + 
    \sum_i | (W_{0})^{(r)}_{i} | 
\leq 
    d_{\mathcal{B}}(W_T) + 
    d_{\mathcal{B}}(W_0) \leq 2 d_{\mathcal{B}}(W_0).
\end{multline*}
Since $\frac{e^{d}}{1 + d}\geq 1$ for any $d>-1$ (recall that $W^{(r)}_l = (\sum_j W_{lj}) -1 > -1$ as well as $W^{(c)}_l = (\sum_i W_{il}) -1 > -1$), then
\[
    \prod_{t=1}^T \frac{1}{1+(W_t)^{(k_t)}_{l_t}} \geq \exp \Big( - \sum_{t=1}^T (W_t)^{(k_t)}_{l_t} \Big)
    \geq \exp(-2 d_{\mathcal{B}}(W_{0})).\qedhere
\]
\end{proof}

\section{Positive diagonals of $W_0$ and relation with the $(W_t)_t$ limit}

\begin{defn}\label{defn:glob_diags}
For an matrix $W\in\Rpos$, and a permutation $\sigma$ of $\{1, 2, \ldots, N\}$, the sequence of elements  $W_{1\sigma(1)},  \ldots, W_{N\sigma(N)}$ is called \emph{the diagonal of $W$ corresponding to $\sigma$}. A positive diagonal of $W$ is a diagonal with positive elements, and $(i_0,j_0)$ is said to lie on a positive diagonal if there  there exists a permutation $ \sigma$ such that $j_0 = \sigma(i_0)$ and $\prod_{i=1}^{N} W_{i, \sigma(i)} > 0$.

A nonnegative matrix \( W \) is said to have \emph{support} if it contains at least one positive diagonal. Formally, $W$ has support if there exists a permutation $ \sigma$ such that $\prod_{i=1}^{N} W_{i, \sigma(i)} > 0$.

A nonnegative matrix \( W \) is said to have \emph{total support} if every positive element of \( W \) lies on a positive diagonal. Formally,
if $ W_{ij} > 0 \implies \exists \sigma \text{ such that } j=\sigma(i)\text{ and } W_{i\sigma(i)} > 0 \text{ for all } i$.

The \emph{positive supported part} of a nonnegative matrix $W$ is the matrix $W^{+}\in\Rpos$ where 
    \[
    (W^{+})_{ij} = 
    \begin{cases}
        W_{ij} & \text{if $(i,j)$ lies on a positive diagonal of $W$}
        \\
        0 & \text{otherwise}.
    \end{cases}
    \]
\end{defn}

\begin{rem}\label{rem:posit_part_total_support}
    If $W\in\Rpos$ has support, its positive supported part has total support. Otherwise, it is the null matrix. 
\end{rem}

The following theorem is a version of the celebrated Birkhoﬀ's Theorem for doubly stochastic matrices. Note that permutation matrices are characterized by their unique positive diagonal; formally, given a permutation $\sigma$ of $\{1,\ldots,N\}$, the permutation matrix $P_\sigma$ is defined as
\[
(P_\sigma)_{ij} = \begin{cases}
    1 & \text{if }j=\sigma(i);
    \\
    0 & \text{if }j\neq \sigma(i).
\end{cases}
\]
\begin{Bthm}
    Every doubly stochastic matrix is a convex combination of permutation matrices.  Formally
    \[
    B \in\mathcal{B} \qquad \iff \qquad B = \sum_{m} \lambda_m P_{\sigma_m}, \quad \lambda_m\geq 0, \sum_m\lambda_m = 1.
    \]
\end{Bthm}
Every doubly stochastic matrix has then total support. A significant converse of this result in this context is given in Theorem~\ref{Sthm}.

\begin{lem}\label{lem:CommClosed}
Given a directed weighted graph $(V,E,w)$, if $W \in \Rpos $ has total support, then each communicating class is closed, that means that the graph is partitioned into strongly connected disjoint components.
\end{lem}
\begin{proof}
Let $V_0 \subseteq V$ be a closed class, that means $W_{ij}=0$ whenever $i\in V_0$ and $j\in V \setminus V_0$. Then for any permutation $\sigma$ with $W_{i\sigma(i)}>0$ we necessary have that $\sigma(i)\in V_0$ whenever $i\in V_0$.
\end{proof}

\begin{defn}
    A directed weighted graph $(V,E,w)$ is called \emph{total irreducible} if $W \in \Rpos $ has total support and is irreducible.
\end{defn}

A star digraph $S_N$ is a weighted directed graph $(V,E,w)$, where $V=\{0,\ldots,,N\}$, $E=\cup_{i=1,\ldots,N} \{(0,i)\cup (i,0)\}$, and $w(e)=1$ for any $e\in E$.
\begin{rem}
    The converse of Lemma~\ref{lem:CommClosed} is not true, since the star graph $S_{N}, N>1$ does not have support.
\end{rem}

\subsubsection*{Diagonals for sub-matrices}
The following definition extends the concept of positive diagonal on the entire set of rows and columns to a subset of them.
\begin{defn}\label{defn:K_diagonal}
Let $\mathcal{K} \subseteq (\{r,c\}\times \{1,\ldots,N\})$. Given $W\in\Rpos$, a \emph{$\mathcal{K}$ diagonal for $W$} is any subset of the elements of $W$ such that if $(k,l)\in\mathcal{K}$, then the $l$-th row (if $k=r$, otherwise the $l$-th column) appears in only one element of the diagonal. A \emph{$\mathcal{K}$-permutation matrix} is a matrix in $\mathbb{R}^{\NtimesN}$ with $1$ on a $\mathcal{K}$ diagonal and $0$ otherwise.

The \emph{$\mathcal{K}$-positive supported part} of a nonnegative matrix $W$ is the matrix $W^{\mathcal{K}}\in\Rpos$ where 
    \[
    (W^{\mathcal{K}})_{ij} = 
    \begin{cases}
        W_{ij} & \text{if $(i,j)$ lies on a positive $\mathcal{K}_0$-diagonal of $W$}
        \\
        0 & \text{otherwise}.
    \end{cases}
    \]
\end{defn}
\begin{rem}\label{rem:partialK}
When $\mathcal{K} = (\{r,c\}\times \{1,\ldots,N\})$, then $\mathcal{K}$-diagonals are the usual diagonals given by permutations. Moreover, there is a natural partial order for diagonals: if $\mathcal{K}'\subseteq \mathcal{K}$, then each $\mathcal{K}$-diagonal is a $\mathcal{K}'$-diagonal. 
Hence, each diagonal given in Definition~\ref{defn:glob_diags} is a $\mathcal{K}$-diagonal, for any $\mathcal{K}$, while the converse is not true.
\end{rem}

\begin{lem}[minimal sub-matrices with no positive diagonals]\label{lem:breveK_0}
    Let $W\in\Rpos$, and assume that for a $\mathcal{K} \subseteq (\{r,c\}\times \{1,\ldots,N\})$ there does not exist a $\mathcal{K}$-positive diagonal.
    Then there exists $\breve{\mathcal{K}} \subseteq \mathcal{K}$, $\breve{\mathcal{K}}\neq\varnothing$ such that there does not still exist a $\breve{\mathcal{K}}$-positive diagonal and
\begin{equation}\label{eq:K_0minimal_intern}
\begin{aligned}
    & \text{$(r,l)\in \breve{\mathcal{K}}$ and $(c,j)\not\in \breve{\mathcal{K}}$} 
&& \implies && 
W_{lj} = 0;
\\
    & \text{$(r,i)\not\in \breve{\mathcal{K}}$ and $(c,l) \in \breve{\mathcal{K}}$} 
&& \implies && 
W_{il} = 0.
\end{aligned}    
\end{equation}    
\end{lem}
\begin{proof}
    If there does not exist a $\mathcal{K}$-positive diagonal, then the partial order relationship given in Remark~\ref{rem:partialK} ensures the existence of a (possibly more than one) minimal $\breve{\mathcal{K}} \subseteq \mathcal{K}$ such that there does not exist a $\breve{\mathcal{K}}$-positive diagonal while there exists a $\mathcal{K}'$-positive diagonal, for any $\mathcal{K}'\subsetneq \breve{\mathcal{K}}$. There always exists a $\mathcal{K}'$-positive diagonal for $\mathcal{K}'$ made by only rows (or columns), since there are no null rows in $W\in\Rpos$, and hence $\breve{\mathcal{K}}\neq\varnothing$.
    
To prove the first assertion of \eqref{eq:K_0minimal_intern}, assume by contradiction that $W_{lj_0} > 0$ for $(r,l)\in \breve{\mathcal{K}}$ and $(c,j_0)\not\in \breve{\mathcal{K}}$.
Then the union of a positive diagonal of $\mathcal{K}' = \breve{\mathcal{K}}\setminus (r,l)$ with $(r,j_0)$ is a positive diagonal of $\breve{\mathcal{K}}$, which is absurd. The same proof may be used to prove the second statement of \eqref{eq:K_0minimal_intern}.
\end{proof}

\begin{rem}\label{rem:max-diag}
By definition, each $\mathcal{K}$-diagonal cannot contain more than $2N-2$ elements.
This limit is achieved for the positive $\mathcal{K}$-diagonal for a star graph $S_{N}, N>1$, where $\mathcal{K}$ contains all the rows and columns except the ones related to the star center. It is obvious that each positive $\mathcal{K}$-diagonal may be chosen with at most $N-1$ equal rows or columns.
\end{rem}

We are now ready to extend the Birkhoff's Theorem to sub-matrices.
\begin{thm}[Birkhoff's theorem for sub-matrices]\label{thm:Birkhoff_subm}
    Let $W\in\Rpos$ and $\mathcal{K} \subseteq \{r,c\}\times\{1,\ldots,N\}$ such that
    \(W^{(k)}_l = 0\) whenever $(k,l)\in \mathcal{K} $. Define $W^*\in \Rpos$ as 
    \[
    W^*_{ij} = \begin{cases}
        W_{ij} & \text{if $(r,i)\in\mathcal{K}$ or $(c,j)\in\mathcal{K}$};
        \\
        0 & \text{if $(r,i)\not\in\mathcal{K}$ and $(c,j)\not\in\mathcal{K}$}.
    \end{cases}
    \]
    Then $W^*$ is a convex combination of $\mathcal{K}$-permutation matrices.
\end{thm}

\begin{proof}
    The proof is an extension of the Birkhoff's algorithm for doubly stochastic matrices.

    As in the proof of Birkhoff's Theorem for doubly stochastic matrices, we will show that there exists a $\mathcal{K}$-permutation matrix $P$ for $W^*$ related to a positive $\mathcal{K}$-diagonal. 
    In this way, if $\underline{w}$ is the smallest non-zero element of the $\mathcal{K}$-diagonal, the matrix $W^*_0 = W^* - P\underline{w}$ is such that 
    \[
    \forall(r,l)\in \mathcal{K}, \qquad \sum_j (W^*_0)_{lj} = 1 - \underline{w}, \qquad \qquad 
    \forall(c,l)\in \mathcal{K}, \qquad \sum_i (W^*_0)_{il} = 1 - \underline{w},  
    \]
since every row and column in $\mathcal{K}$ appears only once in every $\mathcal{K}$-diagonal.
Note that $W^*_0$ has at least one more zero element than $W^*$, and hence this algorithm may be repeated with $W^*$ replaced by $W^*_0 / (1 - \underline{w})$ until $W^*_0$ is the null matrix.

\smallskip 

Then, assume by contradiction that there does not exist a positive $\mathcal{K}$-diagonal for $W^*$.
Apply Lemma~\ref{lem:breveK_0} with $W^*$ and $\mathcal{K}$ to obtain $\breve{\mathcal{K}} \subseteq \mathcal{K}$ for which
\begin{equation}\label{eq:K_minimal_bistochBirkoff}
\begin{aligned}
    & \forall (r,l) \in \breve{\mathcal{K}}  
&& \implies && 
\sum_j W^*_{lj} = 1;
\\
    & \forall (c,l) \in \breve{\mathcal{K}}  
&& \implies && 
\sum_i W^*_{il} = 1,
\end{aligned}
\end{equation}
and such that there is no positive $\breve{\mathcal{K}}$-diagonals.
Equation \eqref{eq:K_0minimal_intern} with \eqref{eq:K_minimal_bistochBirkoff} implies that the number of rows and columns in $\breve{\mathcal{K}}$ are the same.
Moreover, again \eqref{eq:K_0minimal_intern} and \eqref{eq:K_minimal_bistochBirkoff} imply that the sub-matrix of $W^*$ made by the rows and columns in $\breve{\mathcal{K}}$ is doubly stochastic. This implies the existence of a positive $\breve{\mathcal{K}}$-diagonal thanks to Birkhoff's Theorem, which is a contradiction.
\end{proof}

The following theorem states the technical main result, and links the convergence of infinite normalization on one row (column) to the existence of a $\mathcal{K}$-positive diagonal.

\begin{thm}\label{thm:extendAsympNormRowsCols}
    Let $(W_t)_t$ be a normalized sequence generated by $(T^{(k_t)}_{l_t})_{t}$, and let 
    \[
    \mathcal{K}_0 = \{(k,l) \in \{r,c\}\times \{1,\ldots,N\}\colon \text{$k_t = k$ and $l_t = l$ infinitely many times}\}.
    \]
Then, there exists a positive $\mathcal{K}_0$-diagonal for $W_0$ if and only if
\[
\forall (k,l) \in \mathcal{K}_0 \qquad \implies \qquad (W_t)^{(k)}_{l} \mathop{\longrightarrow}\limits_{t\to\infty} 0.
\]
\end{thm}
Before proving this theorem, we give two intermediate results. The first one concerns the relation between the existence of a positive diagonal and the lower bound of the quantity in Proposition~\ref{prop:boundNormsGen}.

\begin{prop}\label{prop:liminfWt}
    Let $(W_t)_t$ be a normalized sequence generated by $(T^{(k_t)}_{l_t})_{t}$, and let 
    \[
    \mathcal{K}_0 = \{(k,l) \in \{r,c\}\times \{1,\ldots,N\}\colon \text{$k_t = k$ and $l_t = l$ infinitely many times}\}.
    \]
If there exists a positive $\mathcal{K}_0$-diagonal for $W_0$, then
% \[
%     \prod_{t=1}^T \frac{1}{1+(W_t)^{(k_t)}_{l_t}} \leq  (\overline{w_0} / \underline{w_0})^{2N-2} N^{N-2},
%     \qquad \text{where} \qquad \begin{aligned}
%         & \underline{w_0} = \min_{(i,j)\colon (W_0)_{ij}>0} (W_0)_{ij}
%         \\
%         & \overline{w_0} = \max_{(i,j)}((W_0)_{ij},1)
%     \end{aligned}
% \]
\[
\sup_T \Big( 
\prod_{t=1}^T \frac{1}{1+(W_t)^{(k_t)}_{l_t}} 
\Big)
< \infty,
\]
while if $(i_0,j_0)$ lies on a positive $\mathcal{K}_0$-diagonal for $W_0$, then
%while if $W_0$ has total support, then
\[
% (W_0)_{i_0j_0} > 0 \qquad \implies \qquad
% \underline{w_0} = \min_{(i,j)\colon (W_0)_{ij}>0} (W_0)_{ij} , \qquad
\inf_T (W_T)_{i_0j_0} > 0 .
\]
\end{prop}
\begin{proof}
Define $\underline{w_0} = \min_{(i,j)\colon (W_0)_{ij}>0} (W_0)_{ij}$ and $\overline{w_0} = \max_{(i,j)}((W_0)_{ij},1)$.
    Let $\{(W_0)_{i_mj_m},m=1,\ldots,M\}$ be a $\mathcal{K}_0$-positive diagonal, and let $\mathcal{I}_0 = \cup\{(r,i_m)\cup(c,i_m),m= 1,\ldots,m\}$ be the sets of the indexes of $\mathcal{K}_0$-diagonal with their multiplicity. By definition of $\mathcal{K}_0$, let $t_0\in\mathbb{N}$ such that $(k_t,l_t)\in \mathcal{K}_0$ for any $t\geq t_0$. Then, if $R_0$ and $C_0$ are the set of different rows and columns in $\mathcal{K}_0$, for $T>t_0$
\begin{multline}\label{eq:pos_K-diag}
\prod_{m=1}^M (W_T)_{i_mj_m} 
= 
\prod_{m=1}^M \Big( (D^{(r)}_T)_{i_mi_m} (W_0)_{i_mj_m} (D^{(c)}_T)_{j_mj_m} \Big)
\\
\begin{aligned}
& \geq 
\Big( \prod_{i_m=1}^M \prod_{t \leq T\colon k_t = r, l_t = i_m} \frac{1}{1+(W_t)^{(r)}_{i_m}} \Big) \underline{w_0}^M 
\Big( \prod_{j_m=1}^M \prod_{t \leq T\colon k_t = c, l_t = j} \frac{1}{1+(W_t)^{(c)}_{j_m}} \Big)
\\
& = 
\Big( \prod_{i \in R_0} \prod_{t \leq T\colon k_t = r, l_t = i} \Big( \frac{1}{1+(W_t)^{(r)}_{i}} \Big)^{\#\{ (r,i) \in \mathcal{I}_0\} } \Big) \underline{w_0}^M 
\Big( \prod_{j \in C_0} \prod_{t \leq T\colon k_t = c, l_t = j} \Big( \frac{1}{1+(W_t)^{(c)}_{j}} \Big)^{\#\{ (c,j) \in \mathcal{I}_0\} } \Big)
\\
& = \underline{w_0}^M 
\prod_{t=1}^{t_0} \Big( \frac{1}{1+(W_t)^{(k_t)}_{l_t}} \Big)^{\#\{ (k_t,l_t) \in \mathcal{I}_0\} }.
\prod_{t=t_0}^T \frac{1}{1+(W_t)^{(k_t)}_{l_t}} . % \Big(\frac{1}{\hat{c}_{W_{t-1} x_t}} \frac{1}{\hat{r}_{T^{(c)}_{x_t} (W_{t-1}) y_t}} \Big) .
\end{aligned}
\end{multline}
The first part of the proof follows immediately, since $(W_T)_{ij} \leq \overline{w_0}$ for any $(i,j)\in \NtimesN$ by Definition~\ref{defn:RC(W)} and hence, by Remark~\ref{rem:max-diag},
\[
\prod_{t=1}^T \frac{1}{1+(W_t)^{(k_t)}_{l_t}} \leq  (\overline{w_0} / \underline{w_0})^M \prod_{t=1}^{t_0} \Big( {1+(W_t)^{(k_t)}_{l_t}} \Big)^{\#\{ (k_t,l_t) \in \mathcal{I}_0\} -1 } \leq (\overline{w_0} / \underline{w_0})^{2N-2} N^{t_0(N-2)}.
\]

For the second part, let $(W_0)_{i_0j_0} > 0$ %. Since $W_0$ has full support, there exists a 
on a $\mathcal{K}_0$-positive diagonal $\{(W_0)_{i_mj_m},m=1,\ldots,M\}$.
Then by Proposition~\ref{prop:boundNormsGen} and \eqref{eq:pos_K-diag}, we have that
\begin{multline*}
\underline{w_0} (\underline{w_0} / \overline{w_0} )^{M-1} \prod_{t=1}^{t_0} \Big( \frac{1}{1+(W_t)^{(k_t)}_{l_t}} \Big)^{\#\{ (k_t,l_t) \in \mathcal{I}_0\} +1}
\exp(-2 d_{\mathcal{B}}(W_{0})) 
\\
\begin{aligned}
& 
\leq  
(1 / \overline{w_0} )^{M-1} \Big (\underline{w_0}^M 
\prod_{t=1}^{t_0} \Big( \frac{1}{1+(W_t)^{(k_t)}_{l_t}} \Big)^{\#\{ (k_t,l_t) \in \mathcal{I}_0\}+1}
\prod_{t=1}^T \frac{1}{1+(W_t)^{(k_t)}_{l_t}} \Big) 
\\
& 
\leq
(W_T)_{i_0j_0} \Big( \prod_{m\neq i_0} \frac{(W_T)_{i_mj_m}}{\overline{w_0}} \Big) \leq  (W_T)_{i_0j_0}. \qedhere
\end{aligned}
\end{multline*}
\end{proof}

Recall Remark~\ref{defn:K_diagonal}: any diagonal is a $\mathcal{K}$-diagonal. Then a sufficient condition for Proposition~\ref{prop:liminfWt} is the fact that $W_0$ has support. In this case, the length of the diagonal is $N$ and $\#\{ (k,l) \in \mathcal{I}_0\}=1$ for any $(k,l)\in \{r,c\}\times \{1,\ldots,N\}$. Then \eqref{eq:pos_K-diag} may be specified to obtain
\[
\prod_{i=1}^N (W_T)_{i\sigma(i)} 
\geq 
\underline{w_0}^N 
\prod_{t=1}^T \frac{1}{1+(W_t)^{(k_t)}_{l_t}} ,
\]
and hence the following corollary holds.
\begin{cor}\label{cor:liminfWt}
    Let $(W_t)_t$ be a normalized sequence generated by $(T^{(k_t)}_{l_t})_{t}$, and let 
    \[
    \mathcal{K}_0 = \{(k,l) \in \{r,c\}\times \{1,\ldots,N\}\colon \text{$k_t = k$ and $l_t = l$ infinitely many times}\}.
    \]
If $W_0$ has support, then
\[
    \prod_{t=1}^T \frac{1}{1+(W_t)^{(k_t)}_{l_t}} \leq  (\overline{w_0} / \underline{w_0})^{N} ,
    \qquad \text{where} \qquad \begin{aligned}
        & \underline{w_0} = \min_{(i,j)\colon (W_0)_{ij}>0} (W_0)_{ij}
        \\
        & \overline{w_0} = \max_{(i,j)}((W_0)_{ij},1)
    \end{aligned}
\]
while if $(i_0,j_0)$ lies on a positive diagonal, then
%while if $W_0$ has total support, then
\[
% (W_0)_{i_0j_0} > 0 \qquad \implies \qquad
% \underline{w_0} = \min_{(i,j)\colon (W_0)_{ij}>0} (W_0)_{ij} , \qquad
\inf_T (W_T)_{i_0j_0} \geq \underline{w_0} (\underline{w_0} / \overline{w_0} )^{N-1} 
\exp(-2 d_{\mathcal{B}}(W_{0})) .
\]
\end{cor}

The second result before the proof of Theorem~\ref{thm:extendAsympNormRowsCols} is a simple technical lemma, and it will be used to link the results given above to the convergence of the normalizations.

\begin{lem}\label{lem:simpleBound}
Let $(d_t)_t$ be a sequence of real number with $-1 < d_t \leq H < +\infty$, $\sup_T|\sum_{=1}^T d_t | \leq K < + \infty$ and 
% $\limsup_t |d_t| > 0$. 
$\sum_{=1}^T d^2_t = \infty$.
Then
\(
\prod_{1}^T (1+d_t) \to 0
\).
\end{lem}
\begin{proof}
For $ 1 \leq s \leq H$ we have $ \tfrac{1-s}{s} \leq \tfrac{1-s}{H} $, while for $ 0 < s < 1$, $\tfrac{1-s}{s} \geq \tfrac{1-s}{H}$. 
In both cases, for any $d\in(-1,H]$
\[
\ln(1+d) - d = \int_1^{1+d} \Big( \frac{1}{s} - 1 \Big) ds \leq \int_1^{1+d} \Big( \frac{1-s}{H} \Big) ds = -\frac{d^2}{2H} ,
\]
since (in both cases) $\sign(s-1) = \sign(d)$, and hence
$$
\prod_{t=1}^T (1+d_t) \leq \exp \sum_{t=1}^T \big( d_t - \tfrac{d_t^2}{2H} \big) \leq \frac{e^K}{\exp\Big( \frac{\sum_{t=1}^T d_t^2}{2H}\Big) }  
\mathop{\longrightarrow}\limits_{T\to\infty} 0. \qedhere
$$
\end{proof}

\begin{proof}[Proof of Theorem~\ref{thm:extendAsympNormRowsCols}]
Let $d_t = (W_t)^{(k_t)}_{l_t}$. By definition and Proposition~\ref{prop:boundNormsGen}, we have then
\[
-1 < d_t \leq \max_{i,j} \big( |(W_t)^{(r)}_{i}| , |W_t)^{(c)}_{j}|| \big) \leq d_{\mathcal{B}}(W_t)\leq d_{\mathcal{B}}(W_0)
\]
and in addition 
\(
    |
    \sum_{t=1}^T d_t
    |
    \leq 2 d_{\mathcal{B}}(W_{0}).
\)
Finally, by Proposition~\ref{prop:liminfWt},
$\inf_T \prod_{t=1}^{T} (1+d_t) > 0$. 

Then we may apply the counternominal of Lemma~\ref{lem:simpleBound} with $H = d_{\mathcal{B}}(W_0)$ and $K = 2H$ to obtain
\begin{equation}\label{eq:RancCto1}
\sum_t d_t^2 <\infty  \qquad \implies \qquad 
\lim_t d_{t} = 
0.
\end{equation}
Take $(r,i_0) \in \mathcal{K}_0$; there exists a subsequence $(t_n)_n$ such that $k_{t_n}=r \text{ and } l_{t_n}= i_0$.
Assume, by contradiction, the existence of $\delta_0> 0$ such that 
\begin{equation}\label{eq:C+to3delta}
\limsup_t |(W_t)^{(r)}_{i_0}| > 2 \delta_0 .
\end{equation}
By \eqref{eq:RancCto1}, there exists $\xi_0 \in\mathbb{N}$ such that $(W_{t_n})^{(r)}_{i_0} < \delta_0$ when $t_n \geq \xi_0$,
and define recursively the sequences $(\eta_t,\xi_t)_t\in\mathbb{N}\cup\{+\infty\}$ as
% $\eta_0 = \xi_0$ and
\begin{align*}
\eta_t 
& =
\begin{cases}
    \inf\{t \geq \xi_{t-1} \colon |(W_t)^{(r)}_{i_0}| \geq 2 \delta_0 \} 
    & \text{ if }\{t \geq \xi_{t-1} \colon |(W_t)^{(r)}_{i_0}| \geq 2 \delta_0 \} \neq \varnothing;
    \\
    +\infty 
    & \text{ if }\{t \geq \xi_{t-1} \colon |(W_t)^{(r)}_{i_0}| \geq 2 \delta_0 \} = \varnothing;
\end{cases}
\\
\xi_t 
& =
\begin{cases}
    \inf\{t \geq \eta_{t} \colon |(W_t)^{(r)}_{i_0}| \leq \delta_0 \} 
    & \text{ if }\{t \geq \eta_{t} \colon |(W_t)^{(r)}_{i_0}| \leq \delta_0 \} \neq \varnothing;
    \\
    +\infty 
    & \text{ if }\{t \geq \eta_{t} \colon |(W_t)^{(r)}_{i_0}| \leq \delta_0 \}  = \varnothing.
\end{cases}
\end{align*}
By induction, if $\xi_{t-1}<\infty$, then \eqref{eq:C+to3delta} implies $\eta_t < \infty$ 
while \eqref{eq:RancCto1} implies $\xi_t < \infty $. 
Moreover, the definition of $\xi_0$ implies that
\begin{equation}\label{eq:notR}
(k_t,l_t) \neq (r,i_0), \qquad \text{ for any $\eta_t < t \leq \xi_t$}.
\end{equation}
Now by \eqref{eq:DeltaW^k_l} and \eqref{eq:notR}, when $\eta_t < t \leq \xi_t$ and $k_t= r$, 
\[
\Delta_{t-1,t} | (W_{\cdot})^{(r)}_{i_0} | - | \Delta_{t-1,t} (W_{\cdot})^{(r)}_{i_0} | = 0.
\]
Hence, since $|a|-|b| - |a-b| \leq 0$, by \eqref{eq:deltaD_B}, 
\begin{multline}\label{eq:B_i_neta}
\Delta_{\eta_t,\xi_t} d_{\mathcal{B}}(W_\cdot) \leq 
\sum_{\eta_t < t \leq \xi_t} 
\Big( 
| (W_{t})^{(r)}_{i_0} | - | (W_{t-1})^{(r)}_{i_0} | -  | (W_{t-1})^{(r)}_{i_0} - (W_{t})^{(r)}_{i_0} |
\Big) 
\\
= | (W_{\xi_t})^{(r)}_{i_0} | - | (W_{\eta_t})^{(r)}_{i_0} |
- \sum_{\eta_t < t < \xi_t} 
 | (W_{t-1})^{(r)}_{i_0} - (W_{t})^{(r)}_{i_0} |
\leq 
\delta_0 - 2 \delta_0
= - \delta_0.
\end{multline}
By Proposition~\ref{prop:boundNormsGen}, $d_{\mathcal{B}}(W_{\xi_0})\leq d_{\mathcal{B}}(W_0)$ and $\Delta_{\xi_{t-1},\eta_t} d_{\mathcal{B}}(W_\cdot) \leq 0$.
Then \eqref{eq:B_i_neta} implies
\begin{align*}
0 \leq d_{\mathcal{B}}(W_{\xi_T}) 
& =
d_{\mathcal{B}} (W_{\xi_0})  +
\bigg[ \sum_{t=1}^T ( d_{\mathcal{B}} W_{\eta_t}) - d_{\mathcal{B}} (W_{\xi_{t-1}}) \bigg]
 +
\bigg( \sum_{t=1}^T d_{\mathcal{B}} (W_{\xi_t})
- d_{\mathcal{B}} (W_{\eta_t}) \bigg)
\\
& \leq d_{\mathcal{B}} (W_{0})  
 + [0]
 - (T \delta_0),
\end{align*}
that is clearly a contradiction for $T$ sufficiently large.

\bigskip

Conversely, if there does not exist a $\mathcal{K}_0$-positive diagonal, then apply Lemma~\ref{lem:breveK_0} with $W_0$ and $\mathcal{K}_0$ to obtain $\breve{\mathcal{K}}_0 \subseteq \mathcal{K}_0$, and assume by contradiction that
\[
\forall (k,l) \in \breve{\mathcal{K}}_0 \qquad \implies \qquad (W_t)^{(k)}_{l} \mathop{\longrightarrow}\limits_{t\to\infty} 0,
\]
that may be read as rows' and columns' sums as
\begin{equation}\label{eq:K_minimal_bistoch}
\begin{aligned}
    & \forall (r,l) \in \breve{\mathcal{K}}_0  
&& \implies && 
\sum_j (W_t)_{lj} \mathop{\longrightarrow}\limits_{t\to\infty} 1;
\\
    & \forall (c,l) \in \breve{\mathcal{K}}_0  
&& \implies && 
\sum_i (W_t)_{il} \mathop{\longrightarrow}\limits_{t\to\infty} 1.
\end{aligned}
\end{equation}
Equation \eqref{eq:K_0minimal_intern} with \eqref{eq:K_minimal_bistoch} implies that the number of rows and columns in $\breve{\mathcal{K}}_0$ are the same.
Moreover, again \eqref{eq:K_0minimal_intern} and \eqref{eq:K_minimal_bistoch} imply that the sequence of square sub-matrices of $(W_t)_t$ made by the rows and columns in $\breve{\mathcal{K}}_0$ converges to $\mathcal{B}$, the set of doubly stochastic matrices. The existence of a converging subsequence implies the existence of a positive $\breve{\mathcal{K}}_0$-diagonal thanks to Birkhoff's Theorem, which is a contradiction.
\end{proof}

The following corollary is Theorem~\ref{thm:extendAsympNormRowsCols} applied to $\mathcal{K} = \{r,c\}\times \{1,\ldots,N\}$.
\begin{cor}\label{cor:convWtToB}
    Let $(W_t)_t$ be a normalized sequence generated by $(T^{(k_t)}_{l_t})_{t}$ for which
for any $i,j\in \{1, \ldots,N\}$, $(k_t,l_t)= (r,i)$ and $(k_t,l_t)= (c,j)$ infinitely many times. Then
\[
\text{$W_0$ has support} \qquad \iff \qquad  d_{\mathcal{B}} (W_{t}) \mathop{\longrightarrow}\limits_{t\to\infty} 0.
\]
\end{cor}
% \begin{rem}
% It is clear from Birkhoff's theorem that no limit to the iteration is possible without support (at least one positive diagonal).
% \end{rem}

\subsection{Characterization of the limit class}
Under the hypothesis of convergence of Corollary~\ref{cor:convWtToB} it is obvious that any converging subsequence of $(W_t)_t$ converges to doubly stochastic matrix, it remains open the question whether the limit exists, and it is unique as function of the sole $W_0$. 
Before stating the answer to this question in Corollary~\ref{cor:limitClass}, we prove that any converging subsequence of $(W_t)_t$ is diagonally equivalent in the more general setting of Theorem~\ref{thm:extendAsympNormRowsCols}.
\begin{thm}\label{thm:limitClass_K0}
    Under the hypothesis of Theorem~\ref{thm:extendAsympNormRowsCols}, if $W_0$ has a positive $\mathcal{K}_0$-diagonal, then the set limits of the sequence $(W_t)_t$ is made by diagonally equivalent matrices to the positive supported part $W_0^{\mathcal{K}_0}$ of $W_0$.
\end{thm}
The proof of this theorem is based on the extension of Birkhoff's Theorem given in Theorem~\ref{thm:Birkhoff_subm} and some arguments given in \cite{Sinkhorn67}, that we report here.
\begin{lem}[{\cite[Lemma~2]{Sinkhorn67}}] \label{Slem:lem2}
    Let $A \in \Rpos$ with total support and
suppose that $(x_n,y_n)_n$
are sequences of positive vectors in $ \mathbb{R}_{>0}^N$ such
that 
\[
A_{ij}> 0 \qquad \implies \qquad (x_n)_i (y_n)_j \text{ converges}.
\]
Then there exist convergent sequences of positive vectors $(x'_n,y'_n)_t$ in $ \mathbb{R}_{>0}^N$ with positive limits such
that $(x_n)_i (y_n)_j = (x'_n)_i (y'_n)_j$ for any $i,j,n$.
\end{lem}
We are not ready to give the proof of our theorem.
\begin{proof}[Proof of Theorem~\ref{thm:limitClass_K0}]
    Let $(W_{t_n})_n$ be a converging subsequence to $A\in\Rpos$, where by \eqref{eq:matreDtRepr}, $W_{t_n} = D^{(r)}_{t_n} W_0 D^{(c)}_{t_n}$, with diagonal matrices $D^{(r)}_{t_n}$ and $D^{(c)}_{t_n}$.
By Proposition~\ref{prop:liminfWt}, if $(i_0,j_0)$ lies  on a positive $\mathcal{K}_0$-diagonal, then $A_{i_0j_0}>0$. Conversely, if $(i_0,j_0)$ does not lie on any positive $\mathcal{K}_0$-diagonal, then $A_{i_0j_0}=0$ by Theorem~\ref{thm:Birkhoff_subm} and Theorem~\ref{thm:extendAsympNormRowsCols}.
In other words
\begin{equation}\label{eq:limit_total_support_K0}
A_{ij} > 0 \qquad \iff \qquad (W_0^{\mathcal{K}_0})_{ij}>0.
\end{equation}
If we define $Z_n = D^{(r)}_{t_n} W_0^{\mathcal{K}_0} D^{(c)}_{t_n}$, it is obvious that $Z_n \to A$. Then \eqref{eq:limit_total_support_K0} implies that
\[
(W_0^{\mathcal{K}_0})_{ij}>0 \qquad \qquad \implies \qquad  (D^{(r)}_{t_n})_{ii} (D^{(c)}_{t_n})_{jj} \mathop{\longrightarrow}\limits_{n\to\infty} 
\frac{A_{ij}}{(W_0^{\mathcal{K}_0})_{ij}},
\]
and hence we may apply Lemma~\ref{Slem:lem2} to $x_n = ((D^{(r)}_{t_n})_{11}, \ldots, (D^{(r)}_{t_n})_{NN})$ and  $y_n = ((D^{(c)}_{t_n})_{11}, \ldots, (D^{(c)}_{t_n})_{NN})$ to obtain that 
\[
A = \widetilde{D}^{(r)} W_0^{\mathcal{K}_0} \widetilde{D}^{(c)}
\]
where $\widetilde{D}^{(r)} = \lim_n \diag(x'_n)$ and $\widetilde{D}^{(c)} = \lim_n \diag(y'_n)$.
\end{proof}

\begin{cor}\label{cor:limitClass}
    Under the hypothesis of convergence of Corollary~\ref{cor:convWtToB}, the sequence $(W_t)_t$ converges to $W_\infty \in \mathcal{B}$ that depends only on the positive supported part $W_0^+$ of $W_0$.
\end{cor}
The proof of this corollary is a consequence of the fact that two doubly stochastic matrices that are diagonally equivalent are equal, as the following result shows. 
\begin{thm}[{\cite[Theorem]{Sinkhorn67}}] \label{Sthm}
    Let $W_0\in\Rpos$. A necessary and
sufficient condition that there exist a doubly stochastic matrix $B$ of
the form $D_1W_0D_2$ where $D_1$ and $D_2$ are diagonal matrices with positive
main diagonals is that $W_0$ has total support. If $B$ exists then it is
unique.
\end{thm}

\section{Applications}

Two key applications of matrix scaling and entropic regularization, within the contexts of directed networks and optimal transport, are presented here. These applications highlight both the theoretical depth and the computational relevance of the paper’s results in modern data science and network analysis.

The first application addresses \textit{Decentralized Random Walks (DRWs)} on time-evolving directed graphs, where edge weights are locally adjusted through normalization at each step. It is shown that under mild conditions, such as irreducibility of the initial support, the transition matrices converge to a doubly stochastic matrix, ensuring uniform visitation frequency across all vertices in the long run.

The second application concerns the approximation of \textit{Optimal Transport (OT)} problems via scaling algorithms under entropic regularization. The OT problem is reformulated to include an entropy term that allows for fast and scalable solutions using the Sinkhorn algorithm. However, our theoretical results highlight a critical limitation: as the regularization parameter \( \varepsilon \to 0 \), numerical underflow can occur due to vanishing entries in the scaled cost matrix \( \xi = e^{-C/\varepsilon} \), resulting in the loss of matrix support and convergence failure.

\subsection{Decentralized Random Walks on Random graphs}
Let $(\Omega,\mathcal{F},(\mathcal{F}_t)_t)$ be a filtered space.

\begin{defn}[Random weighted graphs]
    A \emph{$(\mathcal{F}_t)_t$-adapted directed non-negative weighted graph} (or \emph{random weighted graph}) is a directed graph where each edge has a non-negative adapted process of weight associated with it. Formally, it is a triplet \( \mathcal{G} = (V, E, (W_t)_t) \), where \( G = (V, E) \) is a digraph and for any $t$, $W_t$ is a random weight function on the edges $G$ $\mathcal{F}_t$-measurable.
    % \\
    % A \emph{bistocastic function} is a weight function $w$ on the edges such that
    % \[
    % \forall v_1 \in V, \ \sum_{v \in V \colon (v_1,v)\in E} w(e) = 1, 
    % \qquad 
    % \forall v_2 \in V, \ \sum_{v \in V \colon (v,v_2)\in E} w(e) = 1. 
    % \]
\end{defn}

\begin{defn}[Random Walk on a Random directed Graph and Decentralized Random Walks]
A \emph{random walk on a random directed graph} is a stochastic process $(X_t)_t$ that describes a path consisting of a succession of random steps on the vertices of the random graph $ \mathcal{G} = (V, E, (W_t)_t) $. At each step $t$, the walk moves from the current vertex $X_t$ to a randomly chosen adjacent vertex $X_{t+1}$, following the direction of the edges, and with probability proportional to $W_t((X_t,X_{t+1}))$. 

Given $W_0$ a non-negative $\mathcal{F}_0$-measurable random matrix with support, just before each step $t$, we define
\[
W_t = T^{(r)}_{X_t}(T^{(c)}_{X_t} (W_{t-1}) ),
\]
then the corresponding random walk is called \emph{decentralized random walk (DRW)}.
\end{defn}

% \subsection{DRW}

\begin{prop}\label{prop:linearVisits}
Let $(X_t)_t$ be a DRW on $ \mathcal{G} = (V, E, (W_t)_t) $, and assume that the positive supported part $W_0^+ $ of $W_0\in \Rpos $  is irreducible with probability one.
Then any path $(x_0,\ldots,x_M)$ with $M\geq 0$ and $(W_0^+)_{x_{i-1}x_i}>0$ is visited infinitely often, linearly in $t$.
\end{prop}
\begin{proof}
Since $W_0^+$ is irreducible, for any $i\in \{1,\ldots,N\}$, let $(y_0=i,\ldots, y_{t_0(i)} = x_0)$ a path with $(W_0^+)_{y_{i-1}y_i}>0$ with length less than $N$. 
Define $p_0 = \underline{w_0} (\underline{w_0} / \overline{w_0} )^{N-1} \exp ( - 2 d_{\mathcal{B}} (W_{0})  )$ as in Corollary~\ref{cor:liminfWt}. It is a direct consequence of that corollary that 
\[
P_{W_0^+} (X_t = i, X_{t+1} = y_1, \ldots , X_{t_0(i)} = x_0, X_{t_0(i)+1} = x_1, \ldots, X_{t_0(i)+M} = x_M | X_t = i ) \geq \Big( \frac{p_0}{N} \Big)^{N+M}
\]
holds for any $t= k(N+M)$, $k\in \mathbb{N}$. The thesis follows. 
\end{proof}

\begin{thm}\label{thm:convW_t}
Let $(X_t)_t$ be a DRW on $ \mathcal{G} = (V, E, (W_t)_t) $, and assume $W_0^+ \in \Rpos $  irreducible with probability one.
Then $(W_t)_t$ converges to a double stochastic matrix, that depends only on $W_0^+$.
In addition, 
\[
\forall i\in\{1,\ldots,N\}, \qquad P \Big( \frac{\sum_{t=1}^T \mathbf{1}_{X_t = i}}{T} \mathop{\longrightarrow}\limits_{T\to\infty} \frac{1}{N} \Big) = 1. 
\]
\end{thm}
\begin{proof}
We may apply directly Corollary~\ref{cor:limitClass}, since the hypothesis of Corollary~\ref{cor:convWtToB} are verified by Proposition~\ref{prop:linearVisits}.

In addition, the well-known link between uniform stationary distribution and doubly stochasticity may be used to obtain the last sentence. 
\end{proof}

\subsection{Approximating Transport Problems via Scaling Algorithms and the Role of $\varepsilon$}
\begin{figure}[tbhp]
    \centering
    \includegraphics[width=0.85\linewidth]{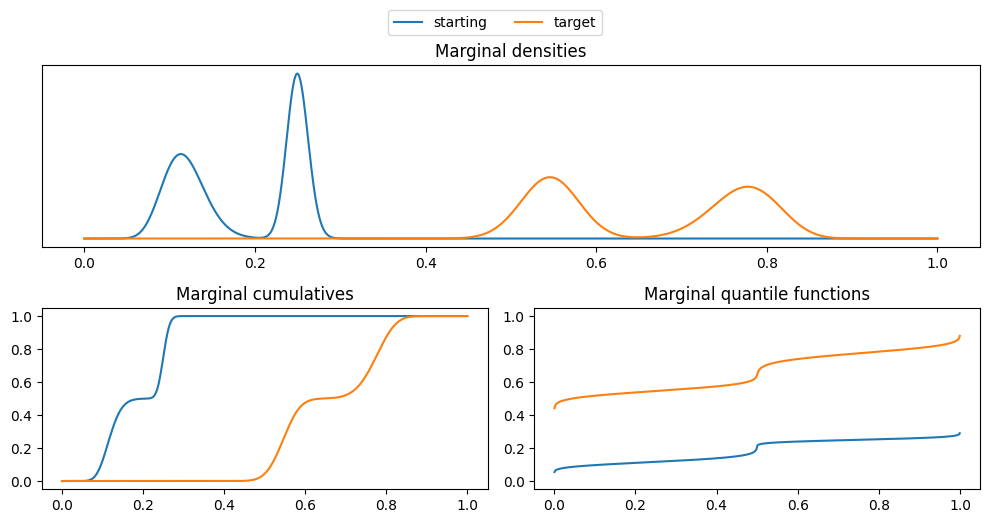}
    \caption{Density (top), cumulative distribution function (bottom-left), and quantile function (bottom-right) of the source and target distributions in the transport problem. In each panel, the blue curve represents the source distribution, and the orange curve represents the target distribution. The source distribution is a mixture of two Beta distributions: \( \text{Beta}(20, 150) \) and \( \text{Beta}(300, 900) \). The target distribution is a mixture of \( \text{Beta}(120, 100) \) and \( \text{Beta}(85, 25) \). Both mixtures use weights \( (1/2, 1/2) \).}
    \label{fig:densities}
\end{figure}
In \cite{Benamou15}, the authors introduce an efficient method to approximate solutions to optimal transport (OT) problems by applying an \emph{entropic regularization} technique. The original linear programming formulation of OT is modified by adding a strongly convex entropy term. This results in a new objective that minimizes a combination of the transport cost and the negative entropy of the coupling matrix. The regularized problem is of the form:
\[
W_\varepsilon(p, q) = \min_{\gamma \in \Pi(p,q)} \langle C, \gamma \rangle - \varepsilon E(\gamma),
\]
where \(\Pi(p, q)\) is the set of couplings with marginals \(p\) and \(q\), \(C\) is the cost matrix, and \(E(\gamma)\) denotes the entropy of \(\gamma\).
A key computational advantage of this regularized formulation is that the optimal coupling \(\gamma^*_\varepsilon\) can be expressed as a diagonal scaling of the matrix \(\xi_{ij} = e^{-C_{ij}/\varepsilon}\). The solution satisfies:
\[
\gamma^*_\varepsilon = \operatorname{diag}(u)\, \xi \, \operatorname{diag}(v),
\]
where the scaling vectors \(u\) and \(v\) are iteratively updated via the Sinkhorn or IPFP algorithm.

However, as \( \varepsilon \to 0 \), the entries of \( \xi = e^{-C/\varepsilon} \) can become exceedingly small, falling below machine precision and resulting in numerical underflow. This causes some entries of \( \xi \) to effectively become zero, leading to a loss of support in the matrix. Consequently, there exists a threshold \( \varepsilon_0 > 0 \) below which \( \xi \) no longer has full support, and no solution is achievable, as stated in Corollary~\ref{cor:limitClass}.

Therefore, while entropic regularization enables fast and stable approximation of transport problems with a modest smoothing parameter \(\varepsilon\), it is not well-suited for approximating the exact (unregularized) OT solution when high accuracy is required.
\begin{figure}[tbhp]
    \centering
    \includegraphics[width=0.85\linewidth]{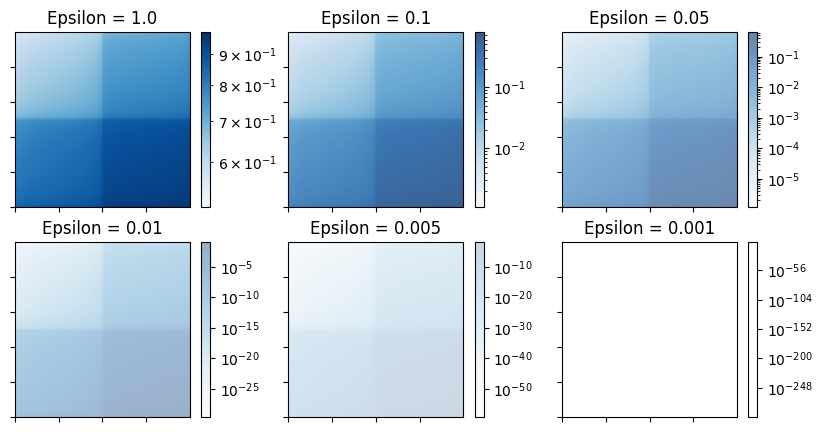}
\caption{Log-scale values of the matrix \( \xi_\varepsilon = e^{-C/\varepsilon} \) for different values of \( \varepsilon \).}
    \label{fig:epsilon}
\end{figure}
Figure~\ref{fig:densities} displays the reconstruction of the two one-dimensional marginal densities (cfr.\ \cite[Fig.~1]{Benamou15}) for the transport problem. The \( N = 1000 \) uniformly spaced quantiles \( (p_n, q_n)_{n=1,\dots,N} \) are used to discretize the continuous densities \( p \) and \( q \). The cost matrix \( C \in \Rpos \) is constructed using the quadratic cost \( C_{ij} = (p_i - q_j)^2 \). The matrix \( \xi_\varepsilon \) is then computed for various values of \( \varepsilon > 0 \). As shown in Figure~2, the positive supported part of \( \xi_\varepsilon \) becomes empty for values of \( \varepsilon \) not too close to zero. For those values of \( \varepsilon \), the sequence of normalizations cannot converge, see Corollary~\ref{cor:limitClass}.

\end{document}